\documentclass[a4paper,10 pt,conference]{ieeeconf}
\usepackage{graphicx}
\usepackage{amssymb}
\usepackage{amsmath}
\usepackage{epstopdf}
\usepackage{algorithmic}
\usepackage{algorithm}

\usepackage{cite}
\usepackage{tikz,pgfplots}
\usepackage{rotating}

\IEEEoverridecommandlockouts 

\newtheorem{theorem}{Theorem}

\title{Approximate Dynamic Programming via \\Sum of Squares Programming}
\author{Tyler H. Summers, Konstantin Kunz, Nikolaos Kariotoglou, \\Maryam Kamgarpour, Sean Summers, and John Lygeros
\thanks{T. H. Summers is supported by an ETH Zurich Postdoctoral Fellowship.}%
\thanks{This research is partially supported by the European Commission under the project MoVeS, FP7-ICT-257005.}
}
\date{\today}                                           

\begin{document}
\maketitle

\begin{abstract}
We describe an approximate dynamic programming method for stochastic control problems on infinite state and input spaces. The optimal value function is approximated by a linear combination of basis functions with coefficients as decision variables. By relaxing the Bellman equation to an inequality, one obtains a linear program in the basis coefficients with an infinite set of constraints. We show that a recently introduced method, which obtains convex quadratic value function approximations, can be extended to higher order polynomial approximations via sum of squares programming techniques. An approximate value function can then be computed offline by solving a semidefinite program, without having to sample the infinite constraint. The policy is evaluated online by solving a polynomial optimization problem, which also turns out to be convex in some cases. We experimentally validate the method on an autonomous helicopter testbed using a 10-dimensional helicopter model. 
\end{abstract}

\section{Introduction}
Many problems in engineering and finance can be modeled as stochastic control problems on infinite state and input spaces, in which a control policy is sought to optimize the behavior of a stochastic dynamical system over a finite or infinite time horizon. While such models are quite general and expressive, the resulting optimization problems are extremely difficult because the decision variable (the control policy) is a \emph{function}, which is generally infinite-dimensional and thus not amenable to computation or even storage on a computer. One general solution method is dynamic programming, which was developed in the seminal work of Bellman in the 1950s \cite{bellman1953bottleneck}. The solution utilizes the Bellman equation, which relates the problem data to the optimal value function and policy. However, solutions of the Bellman equation can be tractably obtained only in a few special cases, when the state and input spaces have very small dimension (hence, can be gridded) or when very strong assumptions are made on the problem data. 

Approximate dynamic programming (ADP) is a collection of heuristic methods for solving stochastic control problems for cases that are intractable with standard dynamic programming methods \cite[Ch. 6]{bertsekas1995dynamic}, \cite{powell2007approximate}. The methods can be classified into three broad categories, all of which involve some kind of function approximation: (1) lookahead/rollout/receding horizon/model predictive control policies, (2) direct policy function approximation, (3) policies based on value function approximation. Here, we will focus on an approach in the last category in which the value function is approximated by a linear combination of pre-specified basis functions. 


A method recently introduced by Wang and Boyd in \cite{wangboyd2010} involves computing an approximate value function by relaxing the Bellman equation to an inequality. One can then obtain a linear program in the basis function coefficients with an infinite set of constraints, one for each state and input pair. In \cite{wangboyd2010}, the value function is approximated with quadratic basis functions, and the S-procedure (see e.g. \cite[Appendix B.2]{boyd2004convex}) is used to transform the infinite constraint set into a linear matrix inequality. The quadratic approximation is used to compute the control policy online, which for certain cases requires solving a small quadratic program. Furthermore, it can be shown that the approximation is a lower bound on the optimal value function and therefore can be used to obtain performance bounds for the system. Finally, the approximation minimizes a norm to the optimal value function, i.e. is in some sense an optimal projection onto the given basis set. The method of Wang and Boyd is based on a similar method from De Farias and Van Roy for finite state and input spaces in \cite{de2003linear}, in which various results on properties of the approximations and constraint sampling methods for finite spaces can be found. Kveton et al also obtain similar results for hybrid spaces \cite{kveton2006solving}.

In this paper, we approximate the value function with more general polynomial basis functions. Polynomial approximations are not new; in fact, Bellman himself studied such approximations in the early 1960s \cite{bellman1963polynomial}. However, the study of polynomial approximations is now particularly interesting in light of recent developments in sum of squares programming \cite{parrilo2003semidefinite} and efficient numerical solvers for semidefinite programming \cite{vandenberghe1996semidefinite,sturm1999using}. In particular, it was shown in \cite{parrilo2003semidefinite} that positivity of a polynomial can be ensured by testing for a sum of squares decomposition, which can then be expressed as a semidefinite program. In stochastic control problems, when the problem data are polynomials, the infinite constraint set in the linear program described above can be expressed as the positivity of a polynomial, which can then be expressed as a semidefinite program. Today, there are widely available numerical solvers for computing solutions to semidefinite programs, making polynomial approximations potentially attractive options for obtaining suboptimal but high-performing solutions to stochastic control problems on high-dimensional spaces.

The main contribution of the present paper is to demonstrate that higher order polynomial value function approximations can be obtained using sum of squares programming techniques. An approximate value function can then be computed offline by solving a semidefinite program, without having to sample the infinite constraint. The policy is evaluated online by solving a polynomial optimization problem; an additional sum of squares constraint can be added to ensure convexity of the approximate value function, making the policy evaluation efficiently computable in certain cases. We also experimentally validate the method using a quadratic value function approximation on an autonomous helicopter testbed. Hover and tracking controllers show good performance, and the policy can be computed online in a few tens of microseconds using recent code generation techniques \cite{mattingley2012cvxgen,Domahidi2012}. To the best of our knowledge, this is the first experimental implementation of this particular approximate dynamic programming method.

The paper is organized as follows. Section II formulates an infinite horizon stochastic control problem, summarizes the dynamic programming solution, and describes an approximate dynamic programming method based on approximating the value function with a linear combination of basis functions. Section III shows how sum of squares techniques can be used to obtain polynomial value function approximations when the problem data are polynomials. Section IV describes an implementation of an approximate dynamic programming controller on an experimental autonomous helicopter testbed. Section V gives concluding remarks and an outlook for future research.

\section{Dynamic Programming and Approximate Dynamic Programming}
In this section, we formulate an infinite horizon stochastic control problem. We then describe the dynamic programming solution and the approximate dynamic programming approach based on value function approximation. 
We use a somewhat informal mathematical style and deliberately set aside non-trivial (but manageable) measurability issues and questions associated with existence of superior non-Markovian or randomized policies to simplify the exposition; see \cite{bertsekas1978stochastic} for a detailed and formal treatment.

\subsection{An Infinite Horizon Stochastic Control Problem}
Consider a discrete-time stochastic dynamical system
\begin{equation}
x^+ = f(x,u,w)
\end{equation}
where $x\in \mathcal{X} \subseteq \mathbf{R}^n$ is the state, $u \in \mathcal{U} \subseteq \mathbf{R}^m$ is the input, and $w\in \mathcal{W} \subseteq \mathbf{R}^p$ is the stochastic process noise with given distribution. We assume that $\mathcal{X}$, $\mathcal{U}$, and $\mathcal{W}$ are closed sets and that $w_t$, $t=0,...$ are independent and identically distributed.

The goal is to find an admissible\footnote{By admissible, we mean that the policy is causal and always satisfies the input constraints. One can formulate problems with state constraints that must be satisfied almost surely or in a probabilistic sense. However, the approximations  described below are not guaranteed to satisfy these; this remains a significant open problem for using these methods. In this paper, state constraints are used to restrict the region of the state space in which the approximation is relevant.} state feedback policy $\psi : \mathcal{X} \rightarrow \mathcal{U}$ with $u_t = \psi(x_t)$ that minimizes the infinite horizon discounted objective function
\begin{equation}
J(x_0,\psi) = \mathbf{E}_w \sum_{t=0}^{\infty} \gamma^t \ell(x_t,u_t)
\end{equation}
where $\gamma \in (0,1)$ is a discount factor and $\ell : \mathcal{X} \times \mathcal{U} \rightarrow \mathbf{R}$ is the stage cost function. This is in general an extremely difficult infinite-dimensional non-convex optimization problem.

The solution can be expressed in principle via dynamic programming. The optimal value function $V^*:  \mathcal{X} \rightarrow \mathbf{R}$ is defined as
\begin{equation}
V^*(x) = \inf_{\psi \in \Psi} J(x,\psi),
\end{equation}
where $\Psi$ is the set of admissible state feedback policies. The function $V^*$ satisfies the Bellman equation
\begin{equation}
 V^*(x) = \inf_{u\in \mathcal{U}} \left\{ \ell(x,u) + \gamma \mathbf{E}_w V^*(f(x,u,w))  \right\}.
\end{equation}
The right-hand-side can be written as an operator on $V^*$
\begin{equation}
V^* = T V^*,
\end{equation}
which has the following properties:
\begin{itemize}
	\item monotonicity: $f \leq g \Rightarrow Tf \leq Tg $
	\item value iteration convergence: $$V^*(x) = \lim_{k\rightarrow \infty} (T^k f)(x) \quad \forall f:\mathcal{X} \rightarrow \mathbf{R}.$$
	\end{itemize}
The optimal policy is then given by
\begin{equation} \label{optpol}
\psi^*(x) = \arg \min_{u\in \mathcal{U}} \left\{ \ell(x,u) + \gamma \mathbf{E}_w V^*(f(x,u,w))  \right\}.
\end{equation}

The optimal value function and policy can only be efficiently computed and stored in a few special cases. When the state, control, and noise spaces are finite and small (approximately $|\mathcal{X}| |\mathcal{U}| |\mathcal{W}| \leq 10^8$), standard dynamic programming methods, e.g. value iteration, policy iteration, or linear programming, can be used to compute the optimal value function and policy. If the state, control, or noise spaces are infinite but have very small dimension (a continuous state space dimension of no more than 4 or 5), the optimal value function and policy can be approximated by gridding the spaces, approximating the dynamics by a finite state Markov chain, and computing the functions at the grid points. For these cases, the dynamic programming methods work for complicated problems but are limited by the so-called curse of dimensionality: computation and storage requirements grow exponentially with the problem dimensions.

For continuous spaces, there is one known case which is computationally tractable in high dimensions. If $\mathcal{X} = \mathbf{R}^n,\ \mathcal{U} = \mathbf{R}^m,\ \mathcal{W} = \mathbf{R}^p$, the dynamics are linear, the noise is Gaussian, and the stage costs are quadratic, then the optimal value function is quadratic, the optimal policy is affine, and both can be computed efficiently from the problem data (either by solving Riccati equations or semidefinite programs). The slightest complication however (e.g. any one of: input/state constraints, non-Gaussian noise, non-quadratic costs, nonlinear dynamics) makes computing the optimal value function extremely difficult. In such cases, we require systematic methods for approximating the optimal value function and policy. 

\subsection{Approximate Dynamic Programming}
We now describe an approximate dynamic programming method that involves approximating the value function with a linear combination of pre-specified basis functions \cite{wangboyd2010,de2003linear}:
\begin{equation}  \label{basisrep}
 \hat{V}(x) = \sum_{i=1}^k \alpha_i \hat{V}_i(x)
\end{equation}
where $\alpha \in \mathcal{A} \subseteq \mathbf{R}^k$. 
The coefficients $\alpha$ are computed offline by solving an optimization problem described below. Then the policy is evaluated online by substituting the approximate value function into the right-hand-side of (\ref{optpol}) (leading to what is sometimes called a greedy policy):
\begin{equation} \label{policy}
\psi^{adp}(x) = \arg \min_{u\in \mathcal{U}} \left\{ \ell(x,u) + \gamma \mathbf{E}_w \hat{V}(f(x,u,w))  \right\}.
\end{equation}
The idea is to find an optimal projection of the optimal value function onto the given set of basis functions and hope that the suboptimal policy (\ref{policy}) yields good performance.

\subsubsection*{Finding the Approximation via Bellman Inequalities}
One method to obtain a value function approximation is to relax the Bellman equation into an inequality \cite{wangboyd2010,de2003linear}. The set of functions that satisfy the Bellman inequality are underestimators of the optimal value function. To see this, suppose a function $\hat{V}$ satisfies $\hat{V} \leq T \hat{V}$. Then by monotonicity of $T$ and value iteration convergence we have
\begin{equation}
\hat{V} \leq T \hat{V} \leq T(T\hat{V}) \leq \cdots \leq \lim_{k\rightarrow\infty} T^k \hat{V} = V^*.
\end{equation}
The Bellman inequality is only a sufficient condition for underestimation, i.e. the set of functions that satisfy the Bellman inequality may not include \emph{all} underestimators. In \cite{wangboyd2010}, it is shown how ``iterated'' Bellman inequalities reduce this conservatism.

This suggests an optimization problem for finding the best value function underestimator in the span of the basis function set:
\begin{equation}
\begin{aligned}
& \text{maximize} \quad \int_\mathcal{X} c(x) \hat{V}(x) dx, \quad c(x) >0 \\
& \text{subject to} \quad \hat{V}(x) \leq T \hat{V}(x), \quad \forall x \in \mathcal{X} 
\end{aligned}
\end{equation}
where $c(x)$ is any positive weighting function. The Bellman inequality constraint is convex in $\alpha$ and can be converted into an infinite set of linear constraints: note that the expression inside the $\inf$ in $T \hat{V}$ is affine in $\alpha$ for each $u$, since expectation is a linear operator, and that the infimum over a family of affine functions is concave and on the right side of the inequality $$ \underbrace{\hat{V}(x)}_\text{affine in $\alpha$} \leq  \underbrace{\inf_{u \in \mathcal{U}} \overbrace{  \left\{ \ell(x,u) + \gamma \mathbf{E}_w \hat{V}(f(x,u,w)) \right\} }^\text{affine in $\alpha$} }_\text{concave in $\alpha$}.  $$ By simply removing the $\inf$ and enforcing the affine constraints for every $u$, we end up with a linear program with an infinite set of constraints, one for every state and input pair:
\begin{equation} \label{semiinflp}
\begin{aligned}
& \text{maximize}  && \int_\mathcal{X} c(x) \hat{V}(x) \\
& \text{subject to}  && \hat{V}(x) \leq \ell(x,u) + \gamma \mathbf{E}_w \hat{V}(f(x,u,w)),\\ & && \forall x \in \mathcal{X}, \ \forall u \in \mathcal{U}\end{aligned}
\end{equation}

The solution to (\ref{semiinflp}) can be shown to be an optimal projection onto the span of the basis functions that satisfy the Bellman inequality in that it minimizes a $c$-weighted $1$-norm to the optimal value function \cite{wangboyd2010,de2003linear}. This gives an interpretation of $c(x)$ as a ``state-relevance weight'', which can be thought of as a distribution over the state space; higher weight can be given to regions of the state space in which we would like better approximation.

In addition to suboptimal policies, value function underestimators also provide performance bounds. For each $x\in \mathcal{X}$ we obtain a number that bounds the performance achievable with \emph{any} admissible state feedback policy. Actual performance of suboptimal policies can be evaluated by Monte Carlo simulation. If the actual performance is close to the bound, the suboptimal policy is close to optimal. Otherwise, either the bound is bad or the policy is bad.

To apply this method in practice, one needs to resolve three main difficulties: (1) evaluating the high-dimensional expectations/integrals, (2) handling the infinite constraint set, and (3) evaluating the policy. Closely related to these issues is the choice of basis functions. We focus on polynomial basis functions. The difficulty in evaluating the expectations/integrals depends on whether or not there are state constraints. When there are no state constraints (or when state constraints are not explicitly taken into account in the integration, as in \cite{wangboyd2010}), all that is needed are the moments of the noise distribution, which we assume to be given or computable from the given distribution. If state constraints are to be explicitly accounted for, certain assumptions must be made on the noise distribution and constraints. If the constraints are boxes and the noise distribution is (approximated by) a polynomial over the constraint set, the integrals can be calculated analytically. The infinite constraint set can be handled in some cases through an appropriate transformation (discussed below); otherwise, one must resort to constraint sampling, incurring more error and losing any theoretical guarantees. To evaluate the policy, an optimization problem (\ref{policy}) must be solved online, which itself is difficult in general even after a good approximate value function is computed. In certain cases (discussed below), this can be made a convex optimization problem and thus amenable to be solved online. 

Wang and Boyd \cite{wangboyd2010} address these issues by using quadratic basis functions. They do not explicitly take into account state constraints in evaluating the integrals, using only the moments of the noise distribution. To handle the infinite constraint, the S-procedure is used to convert it to a linear matrix inequality.  By restricting the quadratic approximation to be convex, the policy can be evaluated online by solving a small quadratic program. In the next section, we show how to extend this approach and approximate value functions by polynomials by utilizing sum or squares programming.





\section{Polynomial Value Function Approximations via Sum of Squares Programming}
In this section, we show how sum of squares programming techniques can be used to find higher-order polynomial value function approximations. Polynomial value functions were studied by Bellman himself in the early 1960s \cite{bellman1963polynomial}. However, recent developments in sum of squares programming and semidefinite programming motivate further study.

\subsection{Sum of Squares and Semidefinite Programming}
In general, determining whether a given multivariate polynomial is everywhere positive is decidable but NP-hard \cite{parrilo2003semidefinite}. An obvious sufficient condition for a polynomial to be positive is for it to be a sum of squares, i.e. be expressible as a finite sum of squared polynomials.
We denote the set of all polynomials in $x\in \mathbf{R}^n$ with real coefficients by $\mathbf{R}[x]$. The set of sum of squares polynomials is denoted by
$$ SOS = \left\{ F(x) \in \mathbf{R}[x] \vert F(x) = \sum_i f_i^2(x), \ f_i(x) \in \mathbf{R}[x] \right\}.$$
and forms a cone in $\mathbf{R}[x]$. 

Let $F(x) \in \mathbf{R}[x]$ be a polynomial of degree $2d$. 
 Let $z$ be a vector of monomials of degree less that or equal to $d$; in particular, $ z(x) = [1,x_1,x_2,...,x_n,x_1x_2,...,x_n^d]^T.$ The following result from \cite{parrilo2003semidefinite} leads to an efficient method to test for a sum of squares decomposition.
\begin{theorem} \label{SOSthm}
The polynomial $F(x)$ is a sum of squares if and only if there exists a symmetric positive semidefinite matrix $Q$ such that $F(x) = z^T(x) Q z(x)$. 
\end{theorem}
By choosing a positive semidefinite matrix subject to affine constraints defined by matching coefficients of $F(x)$ with an expansion of a quadratic form in the monomial vector $z$, one can ensure that $F(x)$ is a sum of squares. This is a semidefinite programming feasibility problem in primal form, which is a convex optimization problem and can be solved efficiently using various well-developed software packages. This method can be used if the coefficients of $F(x)$ must be chosen to satisfy some additional affine constraints.

\subsubsection*{Sum of Squares S-procedure \cite{parrilo2003semidefinite}} 
Now suppose that we want to choose the coefficients of a polynomial so that it is positive over a given semialgebraic set defined by polynomial inequalities, i.e. choose coefficients of $p(x)$ such that $p(x) \geq 0, \ \forall x: g(x) \geq 0$, where $g(x) \in \mathbf{R}[x]$. A sufficient condition for this is the existence of a positive polynomial multiplier $\lambda(x) \geq 0$ such that $p(x) - \lambda(x) g(x) \geq 0$. This can then be restricted to a sum of squares program in which we choose the coefficients of $\lambda$ and $p$ subject to $\lambda(x) \in SOS$ and $p(x) - \lambda(x) g(x) \in SOS$, which is again a semidefinite program. 

\subsection{Polynomial Approximate Dynamic Programming}
Now suppose that all the problem data are polynomial: the dynamics $f$ is polynomial in $x$, $u$, and $w$, the stage cost $\ell$ is polynomial in $x$ and $u$, and the constraint sets $\mathcal{X}$ and $\mathcal{U}$ can be described by polynomial inequalities, i.e. $x \in \mathcal{X}, u \in \mathcal{U} \Leftrightarrow g(x,u) \geq 0$. Suppose also that the value function approximation is a $d$ degree polynomial in $x$ with coefficients as decision variables. Note that none of these polynomials need to be convex. 

The optimization problem (\ref{semiinflp}) can be written
\begin{equation} 
\begin{aligned}
& \text{maximize}  && b^T \alpha \\
& \text{subject to}  && p(x,u,\alpha) \geq 0, \quad \forall g(x,u) \geq 0 \end{aligned}
\end{equation}
where $p(x,u,\alpha)$ is a polynomial in $x$ and $u$ and affine in $\alpha$. Here, we have assumed that the expectations/integrals are either computed analytically or computed based on the moments of the weight and distribution over the state and disturbance space; these computations enter in the entries of $b$ and in the coefficients of $p$. 

We have now a linear optimization problem subject to a constraint that a polynomial is positive on a semialgebraic set, in which all decision variables appear affinely. We can restrict the positivity constraint to a sum of squares constraint using the sum of squares S-procedure to obtain
\begin{equation}  \label{sossdp}
\begin{aligned}
& \text{maximize}  && b^T \alpha \\
& \text{subject to}  && p(x,u,\alpha) - \lambda(x,u) g(x,u) \in SOS \\
&                              && \lambda(x,u) \in SOS.
\end{aligned}
\end{equation}
Via Theorem \ref{SOSthm}, this can be directly transformed into a semidefinite program with variables $\alpha$ and coefficients of $\lambda$, which can be solved offline. 

\subsection{Numerical Example}
To illustrate that higher-order polynomials give improved approximations, we compute quadratic and quartic value function approximations for a 1D problem with linear dynamics, quadratic stage cost, unit-variance Gaussian noise distribution, and state and input constraints:
$$x^+ = x - 0.5 u + w, \quad \ell(x,u) = x^2+u^2, \quad \gamma = 0.99$$ $$\mathcal{X} = [-20,20], \ \mathcal{U} = [-1,1] \Leftrightarrow x^2 \leq 400, \ u^2 \leq 1.$$
Figure 1 shows families of quadratic (blue) and quartic (green) approximations for various weighting functions $c(x)$. The families give similar approximations close to the origin, but the quartic approximations are clearly better near the state constraint boundaries, since all approximations are guaranteed to be underestimators of the optimal value function.

\begin{figure} \label{quadquart}
\resizebox{1\linewidth}{!}{\includegraphics{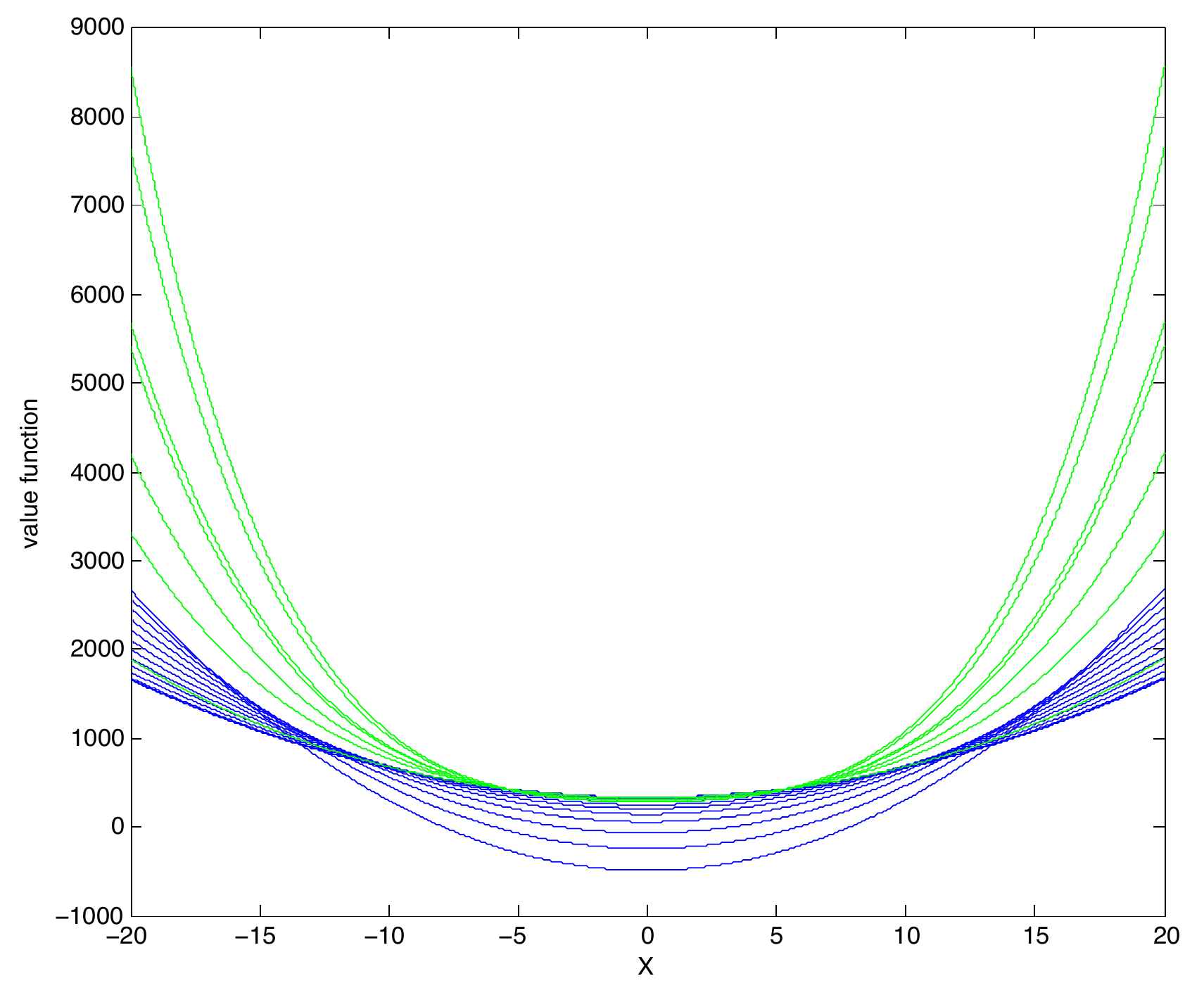}}
\caption{Families of quadratic (blue) and quartic (green) value function approximations.}
\end{figure}

\subsection{Evaluating Policies}
Once we have computed the coefficients of the approximate value function via (\ref{sossdp}), the policy is evaluated online by solving (\ref{policy})
For a given state $x$, this is a polynomial optimization problem in variable $u$ (assuming again the integrals are analytically computable). In general, this is difficult, but we now discuss tractable cases and general techniques.

The policy evaluation (\ref{policy}) is a convex optimization problem whenever $\mathcal{U}$ is convex and the right-hand-side of (\ref{policy}) is a convex function of $u$, which is true if $\hat{V}$ is convex, $\ell$ is convex in $u$, and the dynamics are input affine, i.e. $f(x,u,w) = f_1(x,w) + f_2(x,w)u$, since the composition of a convex function with an affine function is convex, and the sum of convex functions is convex. It turns out that we can enforce convexity of the value function approximation by adding an additional sum of squares constraint on the Hessian of $\hat{V}$ \cite{magnani2005tractable,ahmadi2011complete}. A polynomial $\hat{V}(x)$ is convex if and only if $y^T \nabla^2 \hat{V}(x) y \geq 0, \ \forall x,y \in \mathcal{X}$. Since this is another polynomial inequality, we can ensure convexity of $\hat{V}$ simply by adding the sum of squares constraint $y^T \nabla^2 \hat{V}(x) y \in SOS$ to the semidefinite program (\ref{sossdp}). We have a tradeoff: adding another constraint can only make the approximation worse, even if the optimal value function is convex. However, even if the optimal value function is not known to be convex, it may be advantageous to make the policy evaluation (\ref{policy}) tractable. In these cases, the policy can be evaluated with gradient or interior point methods. Furthermore, since (\ref{policy}) is relatively small (the number of decision variables is the dimension of the input), recent code generation techniques \cite{mattingley2012cvxgen,Domahidi2012} can be used to make the online policy evaluation extremely fast (on the order of microseconds with cheap, embedded processors).


For cases in which a convex approximation is undesirable, there are other techniques that can be considered from a growing literature on general polynomial optimization \cite{lasserre2001global,parrilo2003minimizing,nie2006global}. For example, sum of squares relaxations have been recently explored and observed to recover globally optimal solutions, even in non-convex cases \cite{lasserre2006convergent,parrilo2003minimizing,nie2006minimizing}.

\section{Experimental Validation}
We experimentally validated the approximate dynamic programming method on an autonomous helicopter testbed, which we now briefly describe. Further details can be found in the related paper \cite{kunzetal2013}. 

\subsection{System Architecture}
\begin{figure}
\includegraphics[scale = 0.32]{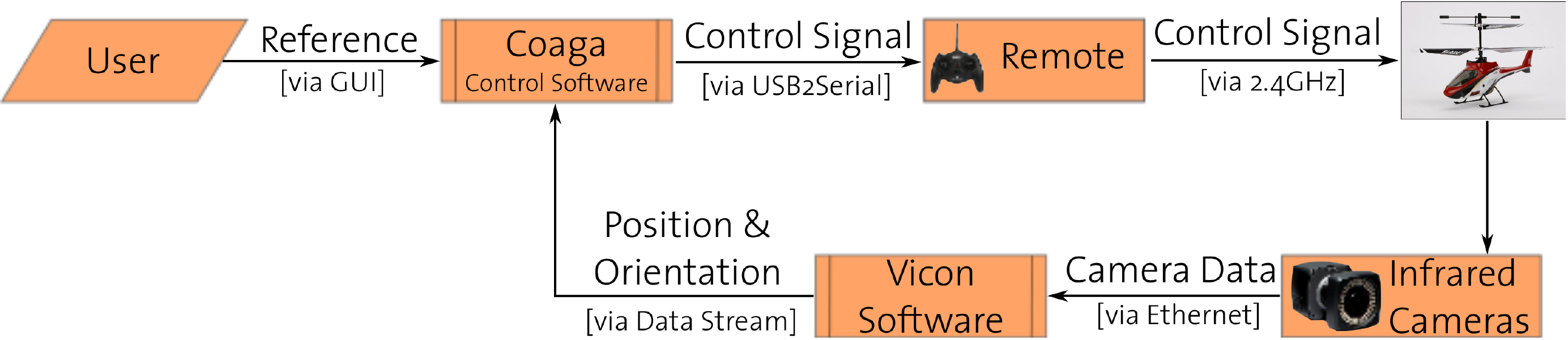}
\caption{The helicopter testbed setup.}
\label{fig:ControlSystem}
\end{figure}

The helicopter navigates in an indoor space equipped with four Vicon Bonita infrared cameras (Vicon product specifications can be found at http://www.vicon.com), which emit infrared light via diodes and record the reflections of markers placed on the helicopter. The camera data is processed by the Vicon Tracker software, which calculates the position and orientation of the helicopter in space. The Coaga \cite{CoagaWiki} custom helicopter control software framework developed by the rCopterX project at ETH Zurich (http://www.rcopterx.ethz.ch) is used to compute the control input from a user-defined desired trajectory and the position and orientation measurements. Coaga is an object-oriented control framework written in C++ to enable fast computations and real-time applications and runs on a stand-alone desktop computer. It also includes a Kalman filter for the estimation of translational and rotational velocities of the helicopter and a feasible trajectory generator based on differential flatness of the model. The computed control input is sent to a remote control via USB. The remote control then transmits the control input directly to the helicopter on a 2.4 GHz frequency. Figure \ref{fig:ControlSystem} shows a block diagram of the experimental setup. 


\begin{figure}
\centering
\includegraphics[scale=0.30]{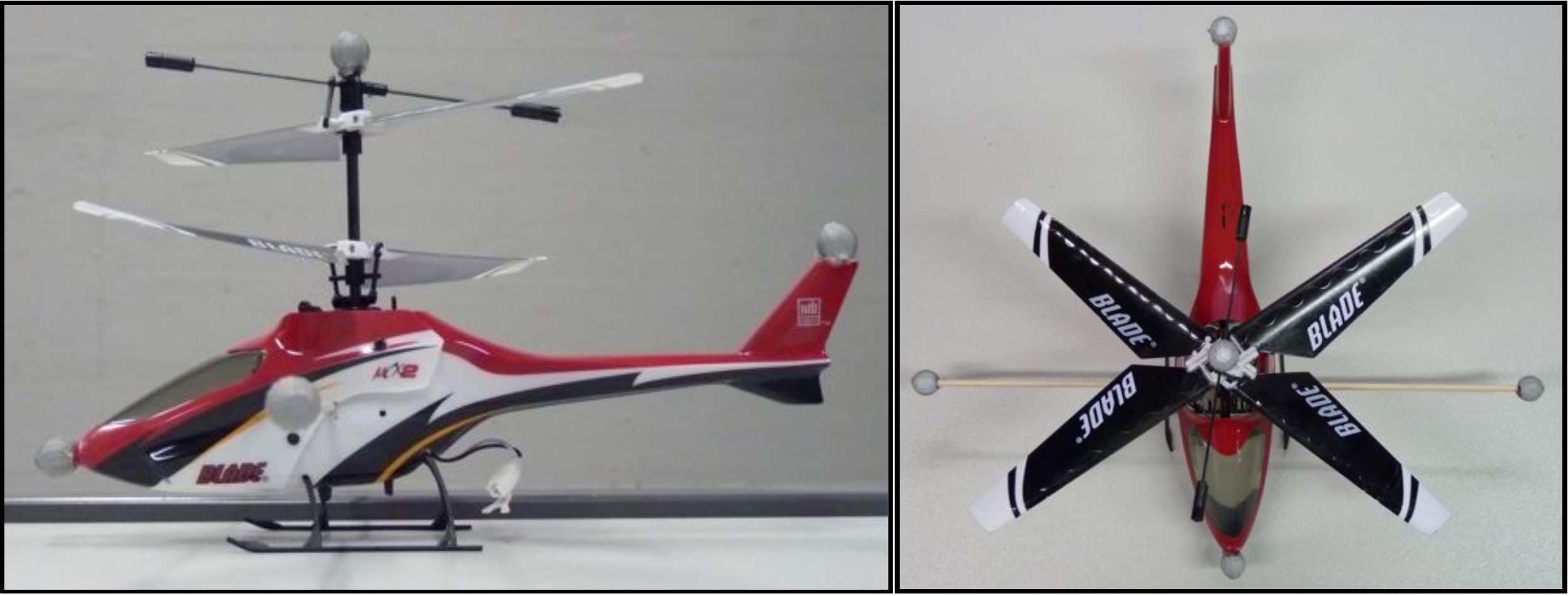}
\caption{Blade mCX2 micro-coaxial helicopters}
\label{fig:BlademCX2}
\end{figure}

The helicopter used in the experiments is a 28-gram Blade mCX2 miniature coaxial helicopter shown in Figure \ref{fig:BlademCX2}. The helicopter is augmented by infrared reflectors so that it can be tracked with the Vicon camera system. It has four control inputs: pitch (for forward flight), roll (for sideways flight), thrust (for vertical flight) and yaw (for heading change). Pitch and roll are actuated with a swashplate mechanism connected to two servos. The thrust is set by the rotation speed of the main rotor motors and yaw is actuated by a rotational speed difference between the two main rotors.

The software used in the experimantal setup consists of Coaga and the Vicon Tracker Software. A convex quadratic value function approximation was computed offline using the semidefinite programming solver SeDuMi \cite{sturm1999using} in the convex optimization modeling framework CVX \cite{grant2008cvx}. We chose a quadratic approximation for validation because this is the simplest type of polynomial that yields a good approximation, and because the resulting quadratic program to be solved online could be done easily with available software and hardware. The ADP online control policy evaluation described above was implemented in C++ and integrated into the Coaga software. The quadratic program was solved with code generated by the CVXGEN code generator \cite{mattingley2012cvxgen} and called in Coaga. 

\subsection{Simplified Helicopter Model}
We used a simplified grey-box model whose parameters were identified based on a least-squares criterion between experimental step-response data and the model prediction. The states used to describe the system are the position and heading in an inertial reference frame $X_{\textit{I}} ,\: Y_{\textit{I}} ,\: Z_{\textit{I}} ,\: \Psi$, the velocity in body-frame $\dot{X}_{\textit{B}} ,\: \dot{Y}_{\textit{B}} ,\: \dot{Z}_{\textit{B}} ,\: \dot{\Psi}$ and two integral states $X_{int} ,\: Y_{int}$ in $x$ and $y$ to reduce steady-state error. The pitch and roll angles of the helicopter are neglected due to the small size of the helicopter. The yaw angle $\Psi$ is defined from an inertial-frame $x$-axis to the body-frame $x$-axis.


The continuous-time state-space representation of the simplified model are
\small
\begin{eqnarray}
\dot{x} &=& Ax+Bu+c \nonumber \\
x &=&
\begin{bmatrix}
X_{\textit{I}} & Y_{\textit{I}} & Z_{\textit{I}} & \Psi & \dot{X}_{\textit{B}} & \dot{Y}_{\textit{B}} & \dot{Z}_{\textit{B}} & \dot{\Psi} & X_{int} & Y_{int}
\end{bmatrix}^T \nonumber \\
u&=&
\begin{bmatrix}
u_x & u_ y & u_z & u_{\Psi}
\end{bmatrix}^T \nonumber \\
c&=&
\begin{bmatrix}
0 & 0 & 0 & 0 & 0 & 0 & -g & 0 & -k_i x_{ref} & -k_i y_{ref}
\end{bmatrix}^T \nonumber
\end{eqnarray}
\footnotesize
\begin{eqnarray} \nonumber
&A
=
\begin{bmatrix}
0 & 0 & 0 & 0 & \cos(\Psi^*) & -\sin(\Psi^*) & 0 & 0 & 0 & 0 \\
0 & 0 & 0 & 0 & \sin(\Psi^*) & \cos(\Psi^*) & 0 & 0 & 0 & 0 \\
0 & 0 & 0 & 0 & 0 & 0 & 1 & 0 & 0 & 0 \\
0 & 0 & 0 & 0 & 0 & 0 & 0 & 1 & 0 & 0 \\
0 & 0 & 0 & 0 & k_x & 0 & 0 & 0 & 0 & 0 \\
0 & 0 & 0 & 0 & 0 & k_y & 0 & 0 & 0 & 0 \\
0 & 0 & 0 & 0 & 0 & 0 & 0 & 0 & 0 & 0 \\
0 & 0 & 0 & 0 & 0 & 0 & 0 & k_{\Psi} & 0 & 0 \\
k_i & 0 & 0 & 0 & 0 & 0 & 0 & 0 & 0 & 0 \\
0 & k_i & 0 & 0 & 0 & 0 & 0 & 0 & 0 & 0 \\
\end{bmatrix} \nonumber \\
&B
=
\begin{bmatrix}
0 & 0 & 0 & 0 \\
0 & 0 & 0 & 0 \\
0 & 0 & 0 & 0 \\
0 & 0 & 0 & 0 \\
b_x & 0 & 0 & 0\\
0 & b_y & 0 & 0 \\
0 & 0 & b_z & 0 \\
0 & 0 & 0 & b_{\Psi} \\
0 & 0 & 0 & 0 \\
0 & 0 & 0 & 0 \\
\end{bmatrix}, \quad \quad
\begin{bmatrix}
b_x \\
k_x \\
b_y \\
k_y \\
b_\Psi \\
k_\Psi \\
b_z \\
\end{bmatrix} \nonumber
=\begin{bmatrix}
2.0 \\
-0.5 \\
2.1 \\
-0.5 \\
111.0 \\
-5.0 \\
18 \\
\end{bmatrix} 
\label{eqa:StateSpace}
\end{eqnarray}
\normalsize
where $\Psi^*$ is the current yaw state around which the system is linearized, the parameters $k_x ,\: k_y \:, k_{\Psi}$ represent fuselage drag, $b_x ,\: b_y \:, b_z ,\: b_{\Psi}$ represent the influence of the inputs $u$ on the states $x$, and $g=9.81 m/s^2$ is gravitational acceleration. These dynamics are time-discretized for implementation. We assumed that unit-variance Gaussian noise acts on accelerations, and we used a quadratic stage cost function.


\subsection{Results}
\begin{figure}
\includegraphics[scale=0.56]{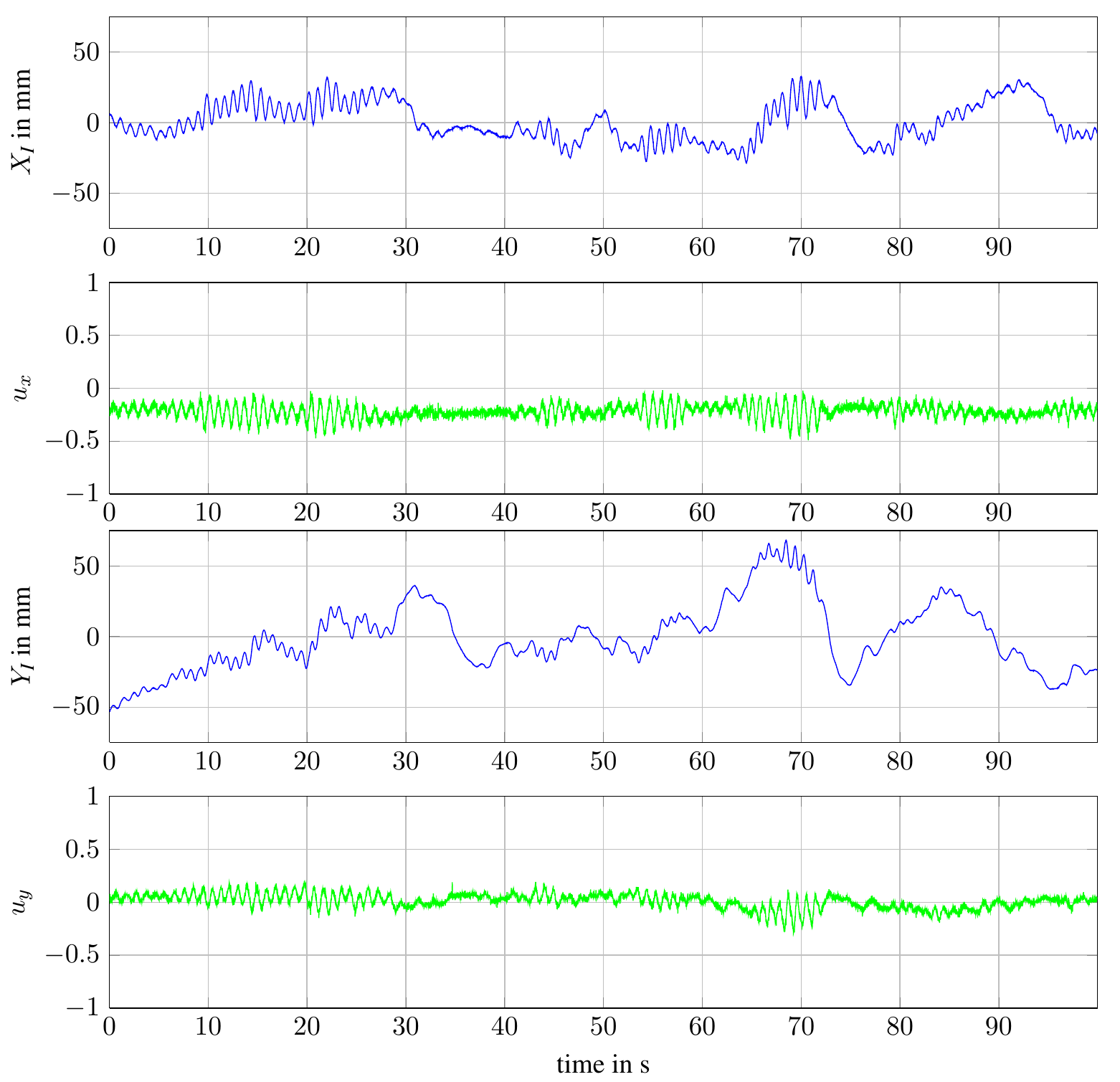}
\caption{ADP controlled helicopter hovering}
\label{fig:HoverFlight}
\end{figure}

\begin{figure}
\includegraphics[scale=0.52]{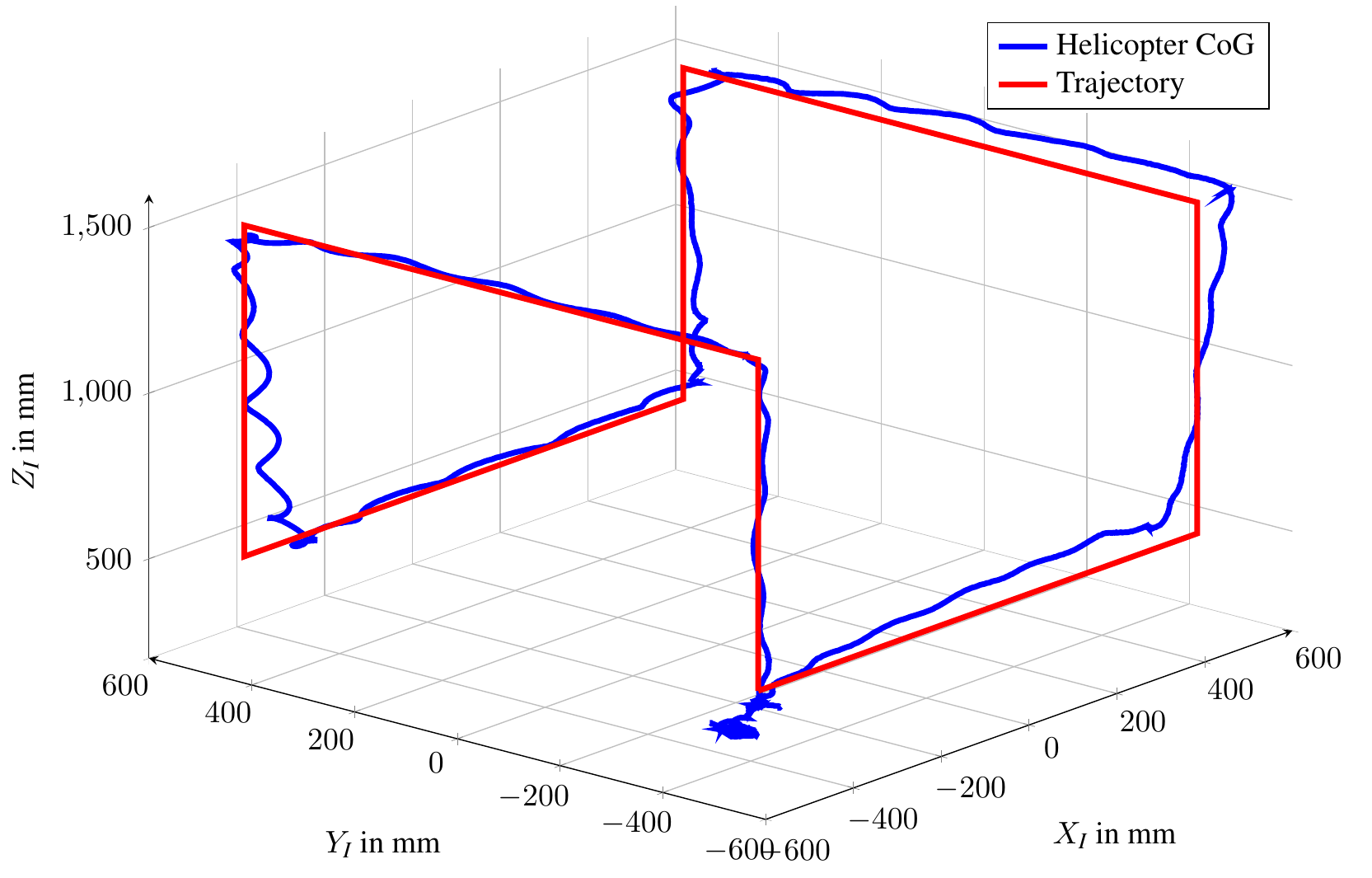}
\caption{ADP controlled helicopter flying a box trajectory}
\label{fig:BoxFlight}
\end{figure}

Both hover and trajectory tracking controllers were developed and tested. The hover performance in $x$ and $y$ can be seen in Figure \ref{fig:HoverFlight}. The maximum deviations in $X_\textit{I}$ and $Y_\textit{I}$ are 3.2cm and 6.8cm, respectively. Trajectory tracking was assessed by having the helicopter fly a box-shaped trajectory with a side length of 1m. A feasible trajectory is precomputed based on the identified helicopter model. The controller tracks both the states and inputs. The results are shown in Figure \ref{fig:BoxFlight}. The results in both hover and trajectory tracking compare favorably with more standard PID and MPC-based controllers reported in \cite{kunzetal2013}. Furthermore, the policy can be evaluated online in a few tens of microseconds, allowing for kilohertz sampling rates (though here we were limited by communication hardware to 50 hertz). 

\section{Concluding Remarks}
We have described an approximate dynamic programming method on infinite state and control spaces. We showed how sum of squares techniques can be used to compute polynomial value function approximations offline via semidefinite programming. The policy is computed online by solving a polynomial optimization problem, which can be made convex in certain cases. Future work will include exploring various application domains, focusing in particular on what can be gained by using higher-order polynomial approximations. Also, the methods can be applied in stochastic reachability problems, which is explored in a companion paper \cite{nikos} via radial basis functions and constraint sampling techniques.


\bibliographystyle{IEEETran}  
\bibliography{refs}  

\end{document}